\documentclass{elsart}

\usepackage{amssymb,amsmath}

\newtheorem{theorem}{Theorem}[section]
\newtheorem{lemma}[theorem]{Lemma}

\newtheorem{proposition}[theorem]{Proposition}
\newtheorem{corollary}[theorem]{Corollary}

\journal{Topology Applications}
\begin{document}

\begin{frontmatter}



\title{The Namioka property of $KC$-functions and Kempisty spaces}

\author{V.V. Mykhaylyuk}

\address{Chernivtsi National University, Department of Mathematical Analysis,
Kotsjubyns'koho 2, Chernivtsi 58012, Ukraine}

\begin{abstract}
A topological space $Y$ is called a Kempisty space if for any
Baire space $X$ every function $f:X\times Y\to\mathbb R$, which is
quasi-continuous in the first variable and continuous in the
second variable has the Namioka property. Properties of compact
Kempisty spaces are studied in this paper. In particular, it is
shown that any Valdivia compact is a Kempisty space and the
cartesian product of an arbitrary family of compact Kempisty
spaces is a Kempisty space.

\end{abstract}

\end{frontmatter}

\section{Introduction}
\label{Introduction} The notion of quasi-continuous mapping, which
was introduced by Kempisty in \cite{K}, occupies an important
place in investigations of the continuity point set of separately
continuous mappings of several variables. Using it,the inductive
pass on the quantity of variables can be reduced to the study of
the continuity point set of two variables mappings which are
quasi-continuous in one variable and continuous in the other
variable.

Let $X$, $Y$ be topological spaces. A function $f:X\times
Y\to\mathbb R$ is called a {\it $KC$-function} if it is
quasi-continuous in the first variable and continuous in the
second variable. The collection of all $KC$-functions $f:X\times
Y\to\mathbb R$ we denote by $KC(X\times Y)$. A function $f:X\times
Y\to\mathbb R$ is called a {\it $\overline{K}C$-function
($\widetilde{K}C$-function)} if it is continuous in the second
variable and there exists a dense set (a dense $G_{\delta}$-set)
$B\subseteq Y$ such that all functions $f_b:X\to \mathbb R$,
$f_b(x)=f(x,b)$, where $b\in B$, are quasi-continuous on $X$. The
collection of all $\overline{K}C$-functions
($\widetilde{K}C$-functions) on $X\times Y$ we denote by
$\overline{K}C(X\times Y)$ ($\widetilde{K}C(X\times Y)$).

A mapping $f:X\times Y\to Z$ has {\it the Namioka property} if
there exists a dense in $X$ $G_{\delta}$-set $A\subseteq X$ such
that $f$ is jointly continuous at each point of $A\times Y$.

A topological space $Y$ is called a {\it Kempisty space} if for
any Baire space $X$ every function $f\in KC(X\times Y)$ has the
Namioka property. This notion was introduced by V.~Maslyuchenko in
\cite{M2}.

 Spaces $Y$ such that for any Baire space $X$ every function
 $f\in \overline{K}C(X\times Y)$
has the Namioka property were studied in \cite{B,M3,N}. The most
general result in this direction was obtained in \cite{M1}. It
gives a possibility to replace the quasi-continuity by the
continuity for those $Y$ all dense subsets of which are separable.
The conditions of countability type (separability, second
countability) on $Y$ are present in all these results.

The Namioka property of $KC$-functions was studied in \cite{m1,m2}
in the case when $X$ is an arbitrary $\alpha$-favorable space. In
particular, it is shown that any $KC$-function on the product of a
$\alpha$-favorable space and a Valdivia compact has the Namioka
property.

On the other hand a compact space $Y$ is called a {\it co-Namioka
space} if for any Baire space $X$ every separately continuous
function $f:X\times Y\to\mathbb R$ has the Namioka property. The
name of these spaces is connected to the classical paper \cite{Na}
and was introduced by G.~Debs in \cite{D}. The class of co-Namioka
compact spaces was intensively studied last time. In particular,
it is shown in \cite{Bo1,Bo2} that any Valdivia compact and the
product of an arbitrary family of co-Namioka spaces are co-Namioka
spaces.

Thus the following questions arise naturally: is any Eberlein,
Corson or Valdivia compact a Kempisty space and is the product of
an arbitrary family of compact Kempisty spaces a Kempisty space?

First, using a dependence of a certain number of coordinates
technique, we obtain a result which implies that a Valdivia
compact is a Kempisty space. Thus, an investigation of the product
of an arbitrary collection of compact Kempisty spaces we reduce to
the case of two factors (these results are announced in \cite{M}).
Then, developing an approach to a study of separately continuous
mappings from \cite{Bo1,Bo2} we show that: $\it {(i)}$ the product
of two compact Kempisty spaces is a Kempisty space; $\it{(ii)}$ in
the definition of a compact Kempisty space one can consider
$KC$-functions which take values in metrizable spaces instead of
real-valued functions, as it was done for co-Namioka spaces in
\cite{Bo1}.

\section{Definitions and an auxiliary assertion}

Let $X, Y$ be topological spaces. A mapping $f:X\to Y$ is called
{\it quasi-continuous at $x_0\in X$} if for any neighborhood $U$
of $x_0$ in $X$ and any neighborhood $V$ of $y_0=f(x_0)$ in $Y$
there exists a nonempty open in $X$ set $G\subseteq U$ such that
$f(G)\subseteq V$. A mapping $f:X\to Y$ is called {\it
quasi-continuous} if it is quasi-continuous at any point $x\in X$.
It is easy to verify that a mapping $f:X\to Y$ is quasi-continuous
if and only if $f(G)\subseteq\overline{f(A)}$ for any nonempty
open in $X$ set $G$ and any dense in $G$ set $A\subseteq X$, where
$\overline{B}$ denotes the closure of $B$.

Let $X$ be a topological space, $Y$ be a metric space with a
metric $d$, $f:X\to Y$ be a mapping, $A\subseteq X$ and
$B\subseteq Y$ be nonempty sets. The oscillation
$\sup\limits_{x',x''\in A}d(f(x'),f(x''))$ of $f$ on $A$ we denote
by $\omega_f(A)$, the oscillation $\inf\limits_{U\in{\mathcal
U}}\omega_f(U)$ of $f$ at $x\in X$ we denote by $\omega_f(x)$,
where ${\mathcal U}$ is the system of all neighborhoods of $x$ in
$X$, and the diameter $\sup\limits_{y',y''\in B}d(y',y'')$ of $B$
we denote by ${\rm diam}\, B$.

For a compact space $K$ by $C(K)$ we mean the Banach space of all
continuous functions $x:K\to \mathbb R$ with the norm $\|x\|=
\sup\limits_{t\in K}|x(t)|$.

For a topological space $X=[0,1]^T$ the space $\Sigma=\{x\in X:
|{\rm supp}\,x|\le\aleph_0\}$, where ${\rm supp }\,x=\{t\in T:
x(t)\ne 0\}$, is called {\it the $\Sigma$-subspace of $X$}. A
compact space $Y$ is called a {\it Corson compact} if it is
embeded homeomorphically into the $\Sigma$-subspace of some
$[0,1]^T$ and a {\it Valdivia compact} if $Y$ is included
homeomorphically into a $[0,1]^T$ such that the $\Sigma$-subspace
of $[0,1]^T$ is dense in the image of $Y$. It is clear that any
Corson compact is a Valdivia compact.

\begin{proposition}
{\it Let $T$ be a set, $(X_t:t\in T)$ be a family of compact
spaces, $X\subseteq \prod\limits_{t\in T}X_t$ be a compact space,
$f:X\to\mathbb R$ be a continuous mapping and $\varepsilon
>0$. Then there exists a finite set $T_0\subseteq T$ such that
$|f(x')-f(x'')|<\varepsilon$ for any $x',x''\in X$ with
$x'|_{T_0}=x''|_{T_0}$.}
\end{proposition}

{\bf Proof.} Using the compactness of $X$ we pick a finite cover
$(U_k: 1\le k\le n)$ of $X$ by open basic sets $U_k$ such that
$|f(x')-f(x'')|<\varepsilon$ if $x',x''\in U_k$ for some $k\le n$.
For each $k\le n$ there exists a finite set $T_k\subseteq T$ such
that $x''\in U_k$ if $x''\in X$ and $x'\in U_k$ with
$x'|_{T_k}=x''|_{T_k}$. It is easy to verify that the set
$T_0=\bigcup\limits_{k=1}^n T_k$ is to be found.$\diamondsuit$

\section{The Namioka property of $KC$-functions}

\begin{theorem}
{\it Let $X$ be a Baire space, $Y\subseteq [0,1]^T$ be a compact,
$\Sigma$ be the $\Sigma$-subspace of a space $[0,1]^T$, $f:X\times
Y\to {\mathbb R}$ be a continuous in the second variable mapping
and $B\subseteq Y\bigcap \Sigma$ be a set with $Y=\overline{B}$
and the function $f_b:X\to {\mathbb R}$, $f_b(x)=f(x,b)$, is
quasi-continuous for any $b\in B$. Then there exists a dense in
$X$ $G_{\delta}$-set $A\subseteq X$ such that $f$ is jointly
continuous at any point of $A\times Y$.}
\end{theorem}

{\bf Proof.} For any $\varepsilon >0$ we put
$A_{\varepsilon}=\{x\in X: \omega_f(x,y)\geq \varepsilon \,\,
\mbox{for some}\,\, y\in Y\}$. Now we prove that all sets
$A_{\varepsilon}$ are nowhere dense in $X$.

Suppose that $\varepsilon >0$ and let $U_0$ be an open in $X$
nonempty set with $U_0\subseteq \overline{A_{\varepsilon}}$.

\begin{lemma}
{\it For any open nonempty set $U\subseteq U_0$ and any set
$S\subseteq T$ with $|S|\leq \aleph_0$ there exist an open
nonempty set $W\subseteq U$ and points $y_1,y_2\in \Sigma\bigcap
Y$ such that $y_1|_S=y_2|_S$ and $|f(x,y_1) - f(x,y_2)|\geq
\frac{\varepsilon}{8}$ for each $x\in W$.}
\end{lemma}

{\bf Proof of Lemma 3.2.} Assume that it is false. Then there
exist an open nonempty set $U\subseteq U_0$ and a set $S\subseteq
T$ with $|S|\leq\aleph_0$ such that for any $y_1, y_2\in
\Sigma\bigcap Y$ with $y_1|_S=y_2|_S$ the set $A(y_1,y_2)=\{x\in
U: |f(x,y_1)- f(x,y_2)|< \frac{\varepsilon}{8}\}$ is dense in $U$.

Consider the continuous mapping $\varphi:Y\to [0,1]^S$,
$\varphi(y)=y|_S$. It is clear that $\varphi(Y)$ is a metrizable
compact. Since $B$ is dense in $Y$ then $\varphi(B)$ is dense in
$\varphi(Y)$. Pick a set $\{b_1, b_2, \dots, b_n, \dots
\}\subseteq B$ such that $\varphi(Y)=\overline{\{\varphi(b_n):n\in
\mathbb N\}}$. Put $\tilde{Y}=\overline{\{b_n:n\in \mathbb N\}}$,
$\tilde{f}=f|_{X\times \tilde{Y}}$. Note that $\tilde{Y}$ is a
metrizable compact, $\tilde{Y}\subseteq \Sigma$ and
$\varphi(\tilde{Y})=\varphi(Y)$. The mapping $\tilde {f}$ belongs
to $\overline{K}C(X,\tilde{Y})$ and has the Namioka property by a
theorem from \cite{B}. Therefore there exists an open nonempty set
$\tilde{U}\subseteq U$ such that ${\rm diam}(\tilde{f}_{\tilde
y}(\tilde{U}))<\frac{\varepsilon}{8}$ for each $\tilde{y}\in
\tilde{Y}$ where $\tilde{f}_{\tilde{y}}:X\to\mathbb R$,
$\tilde{f}_{\tilde{y}}(x)=\tilde{f}(x,\tilde{y})$.

We show that ${\rm
diam}(f_b(\tilde{U}))\leq\frac{3\varepsilon}{8}$ for each $b\in
B$. Fix $b\in B$ and $\tilde{y}\in\tilde{Y}$ so that
$\varphi(b)=\varphi(\tilde{y})$. Put
$\tilde{A}=A(b,\tilde{y})\bigcap \tilde{U}$. Then
$$
|f(a',b)-f(a'',b)|\leq |f(a',b)-f(a',\tilde{y})|+|f(a',\tilde{y})-
f(a'',\tilde{y})|+
$$
$$
+|f(a'',\tilde{y})-f(a'',b)|<\frac{\varepsilon}{8}+\omega_{f_{\tilde{y}}}(\tilde{U})
+\frac{\varepsilon}{8}<\frac{3\varepsilon}{8}
$$
\noindent for any $a', a''\in\tilde{A}$. Thus ${\rm
diam}(f_b(\tilde{A}))\leq\frac{3\varepsilon}{8}$. By the
assumption, $\tilde{A}$ is dense in $\tilde{U}$. Since $f_b$ is
quasi-continuous then
$f_b(\tilde{U})\subseteq\overline{f_b(\tilde{A})}$ and ${\rm
diam}(f_b(\tilde{U}))={\rm
diam}(f_b(\tilde{A}))\leq\frac{3\varepsilon}{8}$ for any $b\in B$.

We prove that $\omega_f(x,y)<\varepsilon$ for any $x\in\tilde{U}$
and $y\in Y$. Fix $x_0\in \tilde{U}$ and $y_0\in Y$. Using the
continuity of $f^{x_0}$, choose a neighborhood $V_0$ of $y_0$ in
$Y$ such that $\omega_{f^{x_0}}(V_0)<\frac{\varepsilon}{16}$. Pick
any points $(x_1,y_1), (x_2,y_2) \in \tilde{U}\times V_0$. Since
$f$ is continuous in variable $y$ and $\overline{B}=Y$ then there
exist points $b_1, b_2\in B\bigcap V_0$ such that
$|f(x_1,y_1)-f(x_1,b_1)|<\frac{\varepsilon}{16}$ and
$|f(x_2,y_2)-f(x_2,b_2)|<\frac{\varepsilon}{16}$. Then
$$ |f(x_1,y_1)-f(x_2,y_2)|\leq
|f(x_1,y_1)-f(x_1,b_1)|+|f(x_1,b_1)-f(x_0,b_1)|+
$$
$$
+|f(x_0,b_1)-f(x_0,b_2)|+|f(x_0,b_2)-f(x_2,b_2)|+|f(x_2,b_2)-f(x_2,y_2)|<
$$
$$
<\frac{\varepsilon}{16}+{\rm
diam}(f_{b_1}(\tilde{U}))+\omega_{f^{x_0}}(V_0) +{\rm
diam}(f_{b_2}(\tilde{U}))+\frac{\varepsilon}{16}\leq
$$
$$
\leq\frac{\varepsilon}{16}+
\frac{3\varepsilon}{8}+\frac{\varepsilon}{16}+\frac{3\varepsilon}{8}+
\frac{\varepsilon}{16}=\frac{15\varepsilon}{16}.
$$
\noindent Therefore $\omega_f(\tilde{U}\times
V_0)\leq\frac{15\varepsilon}{16}$ and
$\omega_f(x_0,y_0)<\varepsilon$. Thus
$A_{\varepsilon}\bigcap\tilde{U}=\O$ but it contradicts
$U_0\subseteq \overline{A_{\varepsilon}}$. Hence, Lemma 3.2 is
proved.$\diamondsuit$

For an arbitrary nonempty open set $U\subseteq U_0$ and set
$S\subseteq T$ with $|S|\leq\aleph_0$ pick a nonempty open set
$\tau(U,S)$ and points $u(U,S), v(U,S) \in \Sigma \bigcap Y$ which
satisfy the conditions of Lemma 3.2.

Describe an strategy $\sigma$ of the second player $\beta$ in the
topological Choquet game on the Baire space $X$. The nonempty open
set $U_0$ is the first move of the player $\beta$. Fix any
nonempty set $S_0\subseteq S$ with $|S_0|\leq \aleph_0$ and for
every nonempty open in $X$ set $V_1\subseteq U_0$ put
$U_1=\sigma(U_0,V_1)=\tau(V_1,S_0)$. Choose any nonempty open in
$X$ set $V_2\subseteq U_1$. Since the points $u_1=u(V_1,S_0)$ and
$v_1=v(V_1,S_0)$ belong to the set $\Sigma$ then
$|T_1|\leq\aleph_0$ where $T_1=({\rm supp}\,u_1)\bigcup ({\rm
supp}\,v_1)$. Denote $S_1=S_0\bigcup T_1$ and put
$U_2=\sigma(U_0,V_1,U_1,V_2)=\tau(V_2,S_1)$.

For any open in $X$ nonempty sets $V_1, V_2,\dots, V_{n-1}$ such
that $U_0\supseteq V_1\supseteq U_1\supseteq \dots \supseteq
U_{n-1}$ where
$U_k=\sigma(U_0,V_1,\dots,U_{k-1},V_{k})=\tau(V_k,S_{k-1})$,
$u_k=u(V_k,S_{k-1})$, $v_k=v(V_k,S_{k-1})$, $T_k=({\rm
supp}\,u_k)\bigcup({\rm supp}\,v_k)$, $S_k=S_{k-1}\bigcup T_k$ by
$k=1,2,\dots,n-1$ and an arbitrary nonempty open in $X$ set
$V_n\subseteq U_{n-1}$ we put
$U_n=\sigma(U_0,V_1,\dots,U_{n-1},V_{n})=\tau(V_n,S_{n-1})$,
$u_n=u(V_n,S_{n-1})$, $v_n=v(V_n,S_{n-1})$, $T_n=({\rm
supp}\,u_n)\bigcup({\rm supp}\,v_n)$, $S_n=S_{n-1}\bigcup T_n$.

Since $X$ is a Baire space then the strategy $\sigma$ is not
winning for the player $\beta$. Therefore there exists a sequence
of nonempty open sets $U_0\supseteq V_1\supseteq U_1\supseteq
\dots \supseteq U_{n-1}\supseteq V_n \supseteq\dots$ where
$U_n=\sigma(U_0,V_1,\dots,U_{n-1},V_{n})$ such that
$\bigcap\limits_{n=0}^{\infty}U_n\ne\O$. Choose a point
$x_0\in\bigcap\limits_{n=0}^{\infty}U_n$. The function $f^{x_0}$
is continuous on the compact $Y\subseteq [0,1]^T$. Therefore by
Proposition 2.1, there exists a finite set $T_0\subseteq T$ such
that $|f(x_0,u)-f(x_0,v)|<\frac{\varepsilon}{8}$ for any $u,v\in
Y$ with $u|_{T_0}=v|_{T_0}$. Since $x_0\in U_n=\tau(V_n,S_{n-1})$,
$u_n=u(V_n,S_{n-1})$, $v_n=v(V_n,S_{n-1})$ then
$|f(x_0,u_n)-f(x_0,v_n)|\geq\frac{\varepsilon}{8}$. Therefore
$u_n|_{T_0}\ne v_n|_{T_0}$ that is for any $n\in \mathbb N$ there
exists a point $t_n\in T_0$ such that $u_n(t_n)\ne v_n(t_n)$. Note
that $t_n\in ({\rm supp}\,u_n)\bigcup({\rm supp}\,v_n)=T_n$, thus
$t_n\in S_{m-1}$ for all $m>n$. Since
$u_m|_{S_{m-1}}=v_m|_{S_{m-1}}$ then $u_m(t_n)=v_m(t_n)$ for all
$m>n$. Hence $t_n\ne t_m$ for all $n\ne m$ and all points $t_k$
are distinct but this contradicts the finiteness of the set
$T_0$.$\diamondsuit$

\begin{corollary}
{\it Let $X$ be a Baire space, $Y$ be a Corson compact and
$f\in\overline{K}C(X\times Y)$. Then $f$ has the Namioka
property.}
\end{corollary}

\begin{corollary}
{\it Let $X$ be a Baire space, $Y$ be a Valdivia compact and
$f\in\widetilde{K}C(X\times Y)$. Then $f$ has the Namioka
property.}
\end{corollary}

{\bf Proof.} Since $Y$ is a Valdivia compact then without loss of
the generality we can assume that $Y\subseteq
\overline{[0,1]^T\bigcap \Sigma}$ where $\Sigma$ is the
$\Sigma$-subspace of $[0,1]^T$. Choose a dense in $Y$
$G_{\delta}$-set $B_1\subseteq Y$ such that all functions
$f_b:X\to \mathbb R$, $b\in B_1$, are quasi-continuous. Put
$B=B_1\bigcap \Sigma$. Since the space $B_2=Y\bigcap\Sigma$ is
countably compact then $B_2$ is a dense Baire subspace of the
Baire space $Y$. Therefore the set $B=B_1\bigcap B_2$ is dense in
$Y$ and by Theorem 3.1, $f$ has the Namioka property.
$\diamondsuit$

\begin{corollary}
{\it Any Valdivia compact is a Kempisty space.}
\end{corollary}

\section{The products of Kempisty spaces}

The following result reduces a study of the products to the case
of two factors.

\begin{theorem} {\it Let $(Y_i: i\in I)$ be a family of compact
spaces $Y_i$ such that for any finite set $I_0\subseteq I$ the
product $\prod\limits_{i\in I_0}Y_i$ is a Kempisty space. Then the
product $Y=\prod\limits_{i\in I}Y_i$ is a Kempisty space.}
\end{theorem}

{\bf Proof.} Consider a Baire space $X$ and $f\in KC(X\times Y)$.
We shall reason similarly as in the proof of Theorem 3.1. Assume
there exist $\varepsilon>0$ and nonempty open set $U_0$ in $X$
such that $U_0\subseteq \overline{A_{\varepsilon}}$ where
$A_{\varepsilon}=\{x\in X: (\exists y\in Y)(\omega_f(x,y)\geq
\varepsilon)\}$. Fix $a\in Y$ and put $\Sigma =\{y\in Y: |\{t\in
T:a(t)\ne y(t)\}|<\aleph_0\}$ (this set will play the same role as
the $\Sigma$-subspace in the proof of Theorem 3.1).

\begin{lemma}
{\it For any open nonempty set $U\subseteq U_0$ and a finite set
$S\subseteq T$ there exist an open nonempty set $W\subseteq U$ and
points $y_1, y_2 \in\Sigma$ such that $y_1|_S=y_2|_S$ and
$|f(x,y_1)-f(x,y_2)|\geq \frac{\varepsilon}{8}$ for every $x\in
W$.}
\end{lemma}

{\bf Proof of Lemma 4.2.} Suppose to the contrary that there exist
an open nonempty set $U\subseteq U_0$ and a finite set $S\subseteq
T$ such that for any $y_1, y_2\in \Sigma$ with $y_1|_S=y_2|_S$ the
set $A(y_1,y_2)=\{x\in U:
|f(x,y_1)-f(x,y_2)|<\frac{\varepsilon}{8}\}$ is dense in $U$.
Denote $\tilde{Y}=\prod\limits_{t\in S}
Y_t\times\prod\limits_{t\in T\setminus S}\{a(t)\}$ and
$\tilde{f}=f|_{X\times\tilde{Y}}$. It is clear that $\tilde{Y}$ is
homeomorphic to $\prod\limits_{t\in S}Y_t$, thus it is a Kempisty
space and $\tilde{f}\in KC(X\times\tilde{Y})$. Hence $\tilde{f}$
has the Namioka property and there exists a nonempty open set
$\tilde{U}\subseteq U$ such that ${\rm
diam}(\tilde{f}_{\tilde{y}}(\tilde{U}))<\frac{\varepsilon}{8}$ for
every $\tilde{y}\in\tilde{Y}$. Then we obtain (analogously as in
the proof of Theorem 3.1) that ${\rm
diam}(f_y(\tilde{U}))<\frac{3\varepsilon}{8}$ for every $y\in
\Sigma$. Now using the continuity of $f$ in the second variable
and the density of $\Sigma$ in $Y$ we obtain that
$\omega_f(x,y)<\varepsilon$ for every $x\in\tilde{U}$ and $y\in
Y$, which is impossible and this completes the proof of the Lemma
4.2. $\diamondsuit$

For any finite set $S\subseteq T$ and a nonempty open set
$U\subseteq U_0$ a nonempty open in $U$ set $W$ and points
$y_1,y_2\in\Sigma$ which satisfy the conditions of Lemma 4.2 we
denote by $\tau(U,S)$, $u(U,S)$ and $v(U,S)$ respectively. We
construct (similarly as in the proof of Theorem 3.1) a winning
strategy for $\beta$ in the Choquet game on the Baire space $X$
(here finite sets $T_k$ are defined by $T_k=\{t\in T: a(t)\ne
u_k(t) \,\,\,\mbox{or}\,\,\, a(t)\ne v_k(t)\}$). It is impossible
because $X$ is a Baire space. $\diamondsuit$

Now we pass to the study of the products of two compact Kempisty
spaces. The following proposition plays a central role in these
investigations. It constitutes a certain development of a somewhat
stronger property of separately continuous functions which has
been used in \cite{Bo2} at similar investigations of co-Namioka
spaces.

Note that for any $\delta >0$ and a metric compact $K$ with a
metric $d$ there exists an integer $m\in \mathbb N$ such that
every set $T\subseteq K$ with $|T|\geq m$ has two distinct
elements $t_1,t_2\in T$ for which $d(t_1,t_2)\leq\delta$ (it is
sufficient to consider a finite $\frac{\delta}{2}$-net in $K$).
The least of such numbers $m$
 we shall call the {\it $\delta$-size of $K$}.

\begin{proposition}
{\it Let $X$ be a Baire space, $Y$, $Z$ be compact spaces,
$P=Y\times Z$, $f:X\times P\to \mathbb R$ be a function which is
continuous in the second variable. Let for any $y\in Y$ the
function $f_y:X\times Z\to \mathbb R$, $f_y(x,z)=f(x,y,z)$, has
the Namioka property. Then for any $\varepsilon >0$ and open in
$X$ nonempty set $U$ there exist functions $b_1, b_2, \dots , b_n
\in C(Z)$ and an open in $X$ nonempty set $U_0\subseteq U$ such
that for every $y\in Y$ there exists a dense in $U_0$ set $A$ such
that for an arbitrary $x\in A$ there is a number $k\in\{1,2,\dots
,n\}$ with $|f(x,y,z)-b_k(z)|\leq\varepsilon$ for each $z\in Z$.}
\end{proposition}

{\bf Proof.} For each $x\in X$ denote by $\varphi_x$ the
continuous mapping $\varphi_x:Y\to C(Z)$,
$\varphi_x(y)(z)=f(x,y,z)$. It is clear that $\varphi_x(Y)$ is a
compact in $C(Z)$.

Fix $\varepsilon > 0$ and an open in $X$ nonempty set $U$. For
every $k\in\mathbb N$ denote by $A_k$ the set of all $x\in U$ such
that $\frac{\varepsilon}{2}$-size of $\varphi_x(Y)$ is not greater
than $k$. Since $\bigcup\limits_{k=1}^{\infty}A_k=U$ and $X$ is a
Baire space then there exist an open in $X$ nonempty set
$U'\subseteq U$ and a number $n$ such that $A_n$ is dense in $U'$
that is $U'\subseteq\overline{A_n}$.

Suppose that for a fixed $\varepsilon$ and an $U$ the conclusion
of Proposition 4.3 is false. In particular, for any open in $X$
nonempty set $U''\subseteq U'$ and a finite set $B\subseteq C(Z)$
there exists an $y\in Y$ such that the set $\bigcup\limits_{b\in
B}\{x\in U'':\|\varphi_x(y)-b\|\leq \varepsilon\}$ is not dense in
$U''$.

We pick an arbitrary point $y_1\in Y$. Since $f_{y_1}$ has the
Namioka property then there exists an open in $X$ nonempty set
$U_1\subseteq U'$ such that
$\|\varphi_{x'}(y_1)-\varphi_{x''}(y_1)\|\leq\frac{\varepsilon}{2}$
for any $x',x''\in U_1$. Fix $x_1\in U_1$ and put
$b_1=\varphi_{x_1}(y_1)$. By the assumption there exists $y_2\in
Y$ such that the set $\{x\in U_1:\|\varphi_{x}(y_2)-b_1\|\leq
\varepsilon \}$ is not dense in $U_1$. Therefore there exists an
open in $X$ nonempty set $U'_1\subseteq U_1$ such that
$\|\varphi_{x}(y_2)-b_1\|>\varepsilon$ for each $x\in U'_1$. Using
the Namioka property of $f_{y_2}$ we find an open in $X$ nonempty
set $U_2\subseteq U'_1$ such that
$\|\varphi_{x'}(y_2)-\varphi_{x''}(y_2)\|\leq
\frac{\varepsilon}{2}$ for any $x',x''\in U_2$.

Given $x_2\in U_2$ we put $b_2=\varphi_{x_2}(y_2)$. Again using
the assumption we find a point $y_3\in Y$ and an open in $X$
nonempty set $U'_2\subseteq U_2$ such that
$\|\varphi_x(y_3)-b_1\|> \varepsilon$ and
$\|\varphi_x(y_3)-b_2\|>\varepsilon$ for every $x\in U'_2$. Now
using the Namioka property of $f_{y_3}$ chose an open in $X$
nonempty set $U_3\subseteq U'_2$ such that
$\|\varphi_{x'}(y_3)-\varphi_{x''}(y_3)\|\leq
\frac{\varepsilon}{2}$ for any $x',x''\in U_3$.

Applying the same arguments $n$ times we obtain a decreasing
finite sequence $(U_k)^{n}_{k=1}$ of open in $X$ nonempty sets
$U_k$ and finite sequences $(x_k)^n_{k=1}$ and $(y_k)^n_{k=1}$ of
points $x_k\in U_k$ and $y_k\in Y$ such that
$\|\varphi_{x'}(y_k)-\varphi_{x''}(y_k)\|\leq
\frac{\varepsilon}{2}$ for any $x',x''\in U_k$ and
$\|\varphi_x(y_k)-b_i\|>\varepsilon$ for arbitrary $x\in U_k$ and
$i\in\{1,2,\dots,k-1\}$ where $b_i=\varphi_{x_i}(y_i)$.

Since $U_n\subseteq U'$ then $A_n\bigcap U_n\ne \O$. Pick an
arbitrary point $x\in A_n\bigcap U_n$ and put $u_k=\varphi_x(y_k)$
for $k=1,2,\dots,n$. Then for $1\leq i<j\leq n$ we have
$$
\|u_j-u_i\| = \|\varphi_{x}(y_j)-\varphi_{x}(y_i)\| \geq
\|\varphi_{x}(y_j)-b_i\| - \|\varphi_{x}(y_i)-\varphi_{x_i}(y_i)\|
> \varepsilon - \frac{\varepsilon}{2} = \frac{\varepsilon}{2}.
$$
\noindent But this contradicts $x\in A_n$.$\diamondsuit$

\begin{theorem}
{\it Let $X$ be a Baire space, $Y$, $Z$ be a compact spaces,
$P=Y\times Z$, $f:X\times P\to \mathbb R$ be a function which is
quasi-continuous in the first variable and continuous in the
second variable and for any $y\in Y$ and $z\in Z$ the functions
$f_y:X\times Z\to\mathbb R$ and $f^z:X\times Y\to \mathbb R$,
$f_y(x,z)=f^z(x,y)=f(x,y,z)$, have the Namioka property. Then $f$
has the Namioka property.}
\end{theorem}

{\bf Proof.} It is sufficient to prove that for any open in $X$
nonempty set $U$ and an $\varepsilon >0$ there exists an open in
$X$ nonempty set $U_0\subseteq U$ such that
$|f(x',y,z)-f(x'',y,z)|<\varepsilon$ for arbitrary $x',x''\in
U_0$, $y\in Y$ and $z\in Z$.

Fix an open in $X$ nonempty set $U$ and $\varepsilon >0$.
Proposition 4.3 implies the existence of functions
$b_1,b_2,\dots,b_n \in C(Z)$ and an open in $X$ nonempty set
$U_1\subseteq U$ such that for any $y\in Y$ there exists a dense
in $U_1$ set $\tilde{A_y}$ such that for every $x\in \tilde{A_y}$
there is a number $k\in \{1,2,\dots,n\}$ with
$|f(x,y,z)-b_k(z)|\leq \frac{\varepsilon}{8}$ for $z\in Z$.

For continuous on $Z$ functions $b_1,b_2,\dots,b_n$ we find a
finite open covering $\mathcal W$ of $Z$ by nonempty sets $W$ such
that $\omega_{b_k}(W)<\frac{\varepsilon}{8}$ for any
$k\in\{1,2,\dots,n\}$ and $W\in \mathcal W$. For every
$W\in\mathcal W$ choose a point $z_W\in W$. Since all functions
$f^{z_W}$ have the Namioka property then there exist an open in
$X$ nonempty set $U_0\subseteq U_1$ such that
$|f(x',y,z_W)-f(x'',y,z_W)|<\frac{\varepsilon}{8}$ for arbitrary
$x',x''\in U_0$, $y\in Y$ and $W\in \mathcal W$. Show that $U_0$
is to be found.

For every $y\in Y$ put $A_y=U_0\bigcap \tilde{A_y}$. It is clear
that $U_0\subseteq \overline{A_y}$ for each $y\in Y$. Fix $y\in Y$
and $W\in \mathcal W$. Recall that for every $a\in A_y$ there
exists a number $k\in \{1,2,\dots,n\}$ such that
$|f(a,y,z)-b_k(z)|\leq \frac{\varepsilon}{8}$ for all $z\in W$.
Since $\omega_{b_k}(W)<\frac{\varepsilon}{8}$ then for arbitrary
$z',z''\in W$ the following inequality holds
$$
|f(a,y,z')-f(a,y,z'')|\leq |f(a,y,z')-b_k(z')| +
|b_k(z')-b_k(z'')| +
$$
$$
+|f(a,y,z'')-b_k(z'')| \leq \frac{\varepsilon}{8} +
\frac{\varepsilon}{8} + \frac{\varepsilon}{8} =
\frac{3\varepsilon}{8}.
$$
\noindent Now for arbitrary $a',a''\in A_y$ and $z\in W$ we have
$$
|f(a',y,z)-f(a'',y,z)|\leq |f(a',y,z)-f(a',y,z_W)| + |f(a',y,z_W)-
$$
$$
- f(a'',y,z_W)| + |f(a'',y,z_W)-f(a'',y,z)| \leq
\frac{3\varepsilon}{8} + \frac{\varepsilon}{8} +
\frac{3\varepsilon}{8} = \frac{7\varepsilon}{8}.
$$

\noindent Since $f$ is quasi-continuous in the first variable and
the set $A_y$ is dense in the open set $U_0$ then
$|f(x',y,z)-f(x'',y,z)|\leq \frac{7\varepsilon}{8} < \varepsilon$
for any $x',x''\in U_0$.$\diamondsuit$

\begin{corollary}
{\it The product of two compact Kempisty spaces is a Kempisty
space.}
\end{corollary}

Theorem 4.1 and Corollary 4.5 imply the following property of
compact Kempisty spaces.

\begin{theorem}
{\it The product of an arbitrary family of compact Kempisty spaces
is a Kempisty space.}
\end{theorem}

\section{$KC$-mappings with values in metrizable spaces}

In this section we carry over the corresponding result from
\cite{Bo1} for separately continuous mappings on the case of
$KC$-mappings which take values in metrizable spaces.

The following statement is an analog of Proposition 4.3.

\begin{proposition}
{\it Let $X$ be a Baire space, $Y$ be a compact space, $Z$ be a
metric space with a metric $d$ and $f:X\times Y\to Z$ be a mapping
which is quasi-continuous in the first variable and continuous in
the second variable. Then for any $\varepsilon >0$ and an open in
$X$ nonempty set $U$ there exist a set $\{z_1, z_2,\dots ,
z_n\}\subseteq Z$ and an open in $X$ nonempty set $U_0\subseteq U$
such that for each $x\in U_0$ and each $y\in Y$ there exists a
number $k\in \{1,2,\dots ,n\}$ such that $d(f(x,y),z_k)\leq
\varepsilon$.}
\end{proposition}

{\bf Proof.} We shall reason similarly like in the proof of
Proposition 4.3.

Fix an open in $X$ nonempty set $U$ and $\varepsilon >0$. Since
$f$ is continuous in the second variable then for any $x\in X$ the
set $Z_x=\{f(x,y): y\in Y\}$ is a metric compact in $Z$. Choose an
open in $X$ nonempty set $U'\subseteq U$, a dense in $U'$ set $A$
and an integer $n\in \mathbb N$ such that
$\frac{\varepsilon}{2}$-size of $Z_a$ is not greater than $n$ for
every $a\in A$.

Suppose that for fixed $\varepsilon$ and $U$ the conclusion of
this proposition  is false. Choose arbitrary points $x_1\in U'$,
$y_1\in Y$ and an open in $X$ nonempty set $U_1\subseteq U'$ such
that $d(f(x,y_1),z_1)<\frac{\varepsilon}{4}$, where
$z_1=f(x_1,y_1)$ for every $x\in U_1$. By the assumption, there
exist points $x_2\in U_1$ and $y_2\in Y$ such that
$d(f(x_2,y_2),z_1)>\varepsilon$. Put $z_2=f(x_2,y_2)$ and using
the quasi-continuity of $f$ in the first variable we find an open
in $X$ nonempty set $U_2\subseteq U_1$ such that
$d(f(x,y_2),z_2)<\frac{\varepsilon}{4}$ for every $x\in U_2$.

Applying the same arguments $n$ times we obtain a decreasing
finite sequence $(U_k)^{n}_{k=1}$ of open in $X$ nonempty sets
$U_k$ and finite sequences $(x_k)^n_{k=1}$ and $(y_k)^n_{k=1}$ of
points $x_k\in U_{k-1}$ for $k=2,\dots, n$ and $x_1\in U'$ and
$y_k\in Y$ such that $d(f(x,y_k),z_k)<\frac{\varepsilon}{4}$ for
every $x\in U_k$ and $d(z_k,z_m)>\varepsilon$ for distinct $k,m\in
\{1,2,\dots,n\}$ where $z_i=f(x_i,y_i)$ for $i=1,2,\dots,n$.

Consider an arbitrary point $a\in A\bigcap U_n$. Then for $1\leq k
< m\leq n$ we have
$$
d(f(a,y_k),f(a,y_m)) \geq d(z_k,z_m) - d(f(a,y_k),z_k) -
$$
$$
- d(f(a,y_m),z_m) > \varepsilon - \frac{\varepsilon}{4} -
\frac{\varepsilon}{4} = \frac{\varepsilon}{2}.
$$
\noindent But this contradicts the choice of $A$.$\diamondsuit$

\begin{theorem}
{\it For any Baire space $X$ and compact space
$Y$ the following conditions are equivalent:

$(i)$ \,\,\,every function $f:X\times Y\to \mathbb R$ which is
quasi-continuous in the first variable and continuous in the
second variable has the Namioka property;

$(ii)$ \,\,\,for any metrizable space $Z$ every mapping $f:X\times
Y\to Z$ which is quasi-continuous in the first variable and
continuous in the second variable has the Namioka property.}
\end{theorem}

{\bf Proof.} It is sufficient to prove the implication
$(i)\Longrightarrow (ii)$. Consider an arbitrary metrizable space
$Z$ and fix a metric $d$ on $Z$ which generates on $Z$ its
topology.

Fix $\varepsilon >0$ and an open in $X$ nonempty set $U$. By
Proposition 5.1 there exist an open in $X$ nonempty set
$U_0\subseteq U$ and a set $\{z_1,z_2,\dots, z_n\}\subseteq Z$
such that for any $x\in U_0$ and $y\in Y$ there exists a number
$k\in \{1,2,\dots, n\}$ with $d(f(x,y),z_k)<
\frac{\varepsilon}{4}$.

For each $k\in \{1,2,\dots, n\}$ the function $f_k:X\times Y\to
\mathbb R$ is defined by $f_k(x,y)=d(f(x,y),z_k)$. By $(i)$ all
functions $f_k$ have the Namioka property. Therefore there exists
an open in $X$ nonempty set $U_1\subseteq U_0$ such that
$|f_k(x',y)-f_k(x'',y)|<\frac{\varepsilon}{4}$ for every $k\in
\{1,2,\dots, n\}$, $x',x''\in U_1$ and $y\in Y$.

Fix $y\in Y$ and $x',x''\in U_1$. Choose $k\in \{1,2,\dots, n\}$
such that $d(f(x',y),z_k)<\frac{\varepsilon}{4}$. Then
$$
d(f(x',y),f(x'',y)) \leq d(f(x',y),z_k) + d(f(x'',y),z_k) =
f_k(x',y) +
$$
$$
+ f_k(x'',y)\leq 2f_k(x',y) + |f_k(x',y)-f_k(x'',y)| <
\frac{3\varepsilon}{4}.
$$
Thus $\omega_f(x,y)<\varepsilon$ for every point $(x,y)\in
U_1\times Y$. Hence $f$ has the Namioka property.$\diamondsuit$



\vspace{3cm}

\begin{center}
V.V. Mykhaylyuk\\
Chernivtsi National University\\
Dept. of Mathematical Analysis\\
Kotsyubyns'koho 2\\
Chernivtsi 58012\\
UKRAINE

\vspace{0.5cm}

mathan@chnu.cv.ua

\vspace{0.5cm}

telephone number: (380) (372) (584888)

\vspace{0.5cm}

\end{center}

\end{document}